\documentclass[letterpaper,12pt]{article}
\usepackage{amsfonts,amsmath,amsthm,graphicx,enumerate,amssymb,amscd,youngtab,tikz,pgf,mathrsfs,pgfplots}
\usetikzlibrary{shapes.geometric, arrows}
\usepackage{pdfpages}
\pgfplotsset{compat=1.15}
\definecolor{uuuuuu}{rgb}{0,0,0}
\definecolor{qqqqff}{rgb}{0,0,1}
\definecolor{ududff}{rgb}{0.30196078431372547,0.30196078431372547,1}
\definecolor{qqttcc}{rgb}{0,0.2,0.8}
\definecolor{ffqqqq}{rgb}{1,0,0}
\definecolor{cqcqcq}{rgb}{0.7529411764705882,0.7529411764705882,0.7529411764705882}
\definecolor{amethyst}{rgb}{0.6, 0.4, 0.8}
\definecolor{amber}{rgb}{1.0, 0.49, 0.0}
 \definecolor{cygreen}{rgb}{0.0, 0.5, 0.0}
\usepackage{array,multirow,makecell}
\setcellgapes{1pt}
\makegapedcells
\newcolumntype{R}[1]{>{\raggedleft\arraybackslash }b{#1}}
\newcolumntype{L}[1]{>{\raggedright\arraybackslash }b{#1}}
\newcolumntype{C}[1]{>{\centering\arraybackslash }b{#1}}
\theoremstyle{definition}

\theoremstyle{remark}

\newcommand{\R}{\mathbb R}

\newcommand{\N}{\mathbb N}

\newcommand{\C}{\mathbb C}

\begin{document}

\title{Accurate algebraic formula for the quintic \\ \& \\Solution by iteration of radicals}
\author{Abdel Missa, Chrif Youssfi}
\maketitle
\newpage
\tableofcontents
\newpage
\begin{center}
{\bf Abstract}\\
\end {center}
According to the Abel-Ruffini theorem [1] and Galois theory [2], there is no solution in finite radicals to the  general quintic equation. This article takes a different approach and proposes a new method to solve the quintic by iteration of radicals. But, the most intriguing result is an accurate algebraic formula for absolute and relative root approximation:\\ $|\textrm{formula- root }|< 4.32\times  \: 10^{-3}$ and $|\textrm{formula/root -1}| < 2.51 \times 10^{-2}$. We then expand some of the geometric properties discussed to construct a trigonometric algorithm that derives all roots. 
\section{Introduction}
Solving polynomial equations occupies a special place in the history of mathematics. The first algebraic solution to the quadratic equation is attributed to Al-Khwarizmi [3] in the 9th century. In the 16th century, a solution to the cubic was found by Del Ferro, Tartaglia and Cardano [4] and was arguably the first significant discovery of the European Renaissance.  Soon thereafter, Ferrari created a method to solve the general quartic equation [5]. The success achieved in finding formulas to lower degree equations stimulated a search to achieve the same ambition for the quintic. Despite efforts by Euler, Gauss, Lagrange and others, no general solutions were found.  Finally, almost 300 years after solving the cubic, Abel provided a complete proof of the impossibility of solving  general quintic polynomials by finite combination of arithmetic operations and radicals. Subsequently, Galois went a step further and provided an exact criterion that characterizes solvable equations:  their Galois group must be solvable.  \\
\\
Naturally, the quintic equation required more tools. The first significant development was the discovery of the Bring-Jerrard normal form  $x^5+ ux+v =0$\\  by using respectively quadratic and quartic Tschirnhaus tranformations [6].  In 1858, Hermite [7] used a simplified normal form to provide the first known solution to general quintic equations in terms of elliptic modular functions.  Shortly afterwards, Brioschi [8] and Kronecker [9] published equivalent solutions. In 1860, Cockle [10] and Harley [11] developed a method for solving the quintic using differential equations, leading to a solution in hypergeometric functions.  And in 1884, Klein [12] provided an icosahedral solution of the quintic.  More recently, Doyle and McMullen [13] solved the quintic using an iterative method, while Glasser [14] provided a derivation to trinomial equations that leads to the same hypergeometric solution for the quintic as Cockle and Harley. \\
\\
In this article, we explore a new method that solves the quintic by iteration of radicals. Perhaps the most astonishing finding is the ability to find an accurate global approximation in radicals to one the five roots from the first iteration: $|\textrm{1st iteration- root} | < 4.32\: 10^{-3}$. We then proceed by proving the speed of convergence of the proposed iterative algorithm. \\
\\
Next, we discuss Bring radicals. Recall that the Bring radical of a complex number $a$ is any of the five roots of: $x^5+x+a=0$. More importantly, quintic equations can be solved in closed form using radicals and Bring radicals. One might think that a simple iterative algorithm $x_{k+1}=-\sqrt[5]{a+x_k}$ can solve the Bring equation. However, this is not the case as we show a counter example (when $a$ is real) where the algorithm does not converge regardless of the starting point's choice outside of the real root itself. On other hand, we prove that our proposed iterative algorithm converges globally for any complex number $a$. Furthermore, we provide a radical approximation to one of the Bring radicals:  $|BR(z) – \textrm{approximation }|  < 2.90  \: 10^{-2}$.
\\
\\
Finally, we expand some of the properties used in the proposed algorithm to provide the location of the five roots in the complex plane. This leads to a trigonometric bisection algorithm that solves general quintic equations.

\section{Reduced Forms}
The general quintic equation can be transformed to the Bring–Jerrard normal form:
\[ v^5+d_1v+d_0=0\]
From hereon, we define the n-th root of a complex number as the one with its argument in $[-\pi/n, \pi/n[$. \\
\\
Excluding the trivial cases of $d_1=0$ and $d_0=0$, a change of variable $x=v/\sqrt[4]{d_1}$ can further simplify the equation to Form 1: 
\begin{equation}
x^5+x+a=0 \qquad \textrm{where}\qquad a=\frac{d_0}{\sqrt[4]{d_1^5}} \qquad (a\neq 0)
\end{equation}
\\
Another change of variable $z=a/x$  leads to Form 2:
\begin{equation}
\frac{z^5+z^4}{2}=\lambda \qquad \textrm{where} \qquad \lambda=-\frac{a^4}{2}
\end{equation}
Let $y=uz$, (2) becomes:
\[\frac{y^5+u y^4}{2}=u^5\lambda \]
Choose $u$ to satisfy: \[u^5 =\frac{|\lambda|}{\lambda} \qquad \textrm{or} \qquad u =\left(\frac{|\lambda|}{\lambda}\right)^{\frac{1}{5}}=e^{i\theta} \qquad \theta \in [-\frac{\pi}{5},\frac{\pi}{5}[.\]
This leads to Form 3:
\begin{equation}
\frac{y^5+uy^4}{2}=\xi \qquad \xi =|\lambda| \in \R^{+*}
\end{equation}
Think of Form 3 as a simple rotation of Form 2 with the benefit of solving the equation close to the real axis.\\
\\
Notice that $s$ is solution of (3) if and only if its conjugate $\bar{s}$ is a solution of 
\[\frac{y^5+\bar{u}y^4}{2}=\xi  \]
Therefore we can assume from hereon that $0\leq \theta \leq \pi/5$.

\section{Iteration of radicals algorithm to solve equations in Form 3}

\subsection{Algorithm}
Let us consider Equation (3) 
\[\frac{y^5+uy^4}{2}=\xi \qquad \xi \in \R^{+*} \qquad (3)\]
\[u=e^{i\theta} \qquad  \theta \in \left[0, \frac{\pi}{5}\right]\]
It can be exploited in two ways. First, by factorizing:
\begin{equation}
y=\sqrt[4]{\frac{2\xi}{u+y}} 
\end{equation}
Second, by completing the quintic:
\begin{equation} 
\left(y+\frac{u}{5}\right)^5=2\xi+\frac{2u^2}{5}y^3 +\frac{2u^3}{25}y^2+\frac{u^4}{125}y+\frac{u^5}{3125} 
\end{equation}
By substituting $y$ in the right-hand side of $(5)$:
\begin{equation} 
y=G(\xi,y)
\end{equation}
with:
\begin{equation} 
G\left(\xi,y\right)=\sqrt[5]{2\xi+\frac{2u^2}{5}t^3 +\frac{2u^3}{25}t^2+\frac{u^4}{125}t+\frac{u^5}{3125}}-\frac{u}{5}
\end{equation}
\begin{equation} 
t=\sqrt[4]{\frac{2\xi}{u+y}}
\end{equation}
This inspires the following fixed point algorithm with the following sequence: $(y_k)_{k \in \N}$ by: \\
\\
\underline {Starting point}\footnote{The intuition behind $\alpha$ is discussed below.}
\begin{equation}
y_0=\left(\frac{\xi}{\alpha}\right)^{\frac{2}{9}} \textrm{  with  }\alpha=\sqrt{\frac{2+\sqrt{2}}{4}} 
\end{equation} 
\underline {Iterative process}
\begin{equation}
y_{k+1}=G\left(\xi,y_k\right)
\end{equation}
\\
To prove the convergence and the accuracy of this algorithm, we proceed with the following steps:
\begin{itemize}
\item Section 3.2 shows that there exists a unique root $y^*$ for Equation (3) near the positive real axis.  
\item Section 3.3  shows that $\forall\: \xi \in \R^{+*}$ and $\forall \:\theta \in \left[0,\pi/5\right]$: \[\left|y_1-y^{*}\right| <4.32 \textrm{ } 10^{-3}\]
In other words, the closed form formula: 
\[y_1=\sqrt[5]{2\xi+\frac{2u^2}{5}t^3 +\frac{2u^3}{25}t^2+\frac{u^4}{125}t+\frac{u^5}{3125}}-\frac{u}{5}\]
with
\[t=\sqrt[4]{\frac{2\xi}{u+\left(\frac{\xi}{\alpha}\right)^{\frac{2}{9}}}}\]

provides an accurate absolute approximation of $y^*$
\item Section 3.4  shows that the relative error of the formula is also small: $\forall\: \xi \in \R^{+}$ and $\forall \:\theta \in \left[0,\pi/5\right]$ \[\left|\frac{y_1}{y^{*}}-1\right| <2.58 \textrm{ } 10^{-2}\]
\item Finally, section 3.5 shows the sequence $(y_k)_{k \in \N}$ converges toward $y^{*}$ and:
\[\forall \: k \in \N \qquad \left|y_{k+1}-y^*\right| < \frac{1}{K}\left|y_{k}-y^*\right| \qquad  K\approx 15.44\]
\end {itemize}

\subsection{Unique root near the positive real axis}
Consider equation (3). When $\theta=0$, $u=1$ and it is obvious that there is a unique real solution that satisfies (3).
Now assume $\theta >0$. We will show in this section that there exists a unique root $y$ of Equation (3) near the real axis: 
\begin{equation} 
y= r e^{i\sigma} \qquad \textrm{with } \qquad r>0  \qquad \textrm{and } \qquad {\bf \sigma \in \left[-\frac{\theta}{4},0\right[}
\end{equation}
The real and the imaginary parts of the equation are:
\begin{equation}
r^5 \cos\left(5\sigma\right)+r^4 \cos\left(\theta+4 \sigma\right)=2\xi
\end{equation}
and
\begin{equation}
r^5 \sin\left(5\sigma\right)+r^4 \sin\left(\theta+4 \sigma\right)=0
\end{equation}
From (13):
\begin{equation}
r=-\frac{\sin\left(\theta+4 \sigma\right)}{\sin\left(5\sigma\right)}
\end{equation}
Substituting $r$ in (12):
\begin{equation}
f(\sigma)=\frac{\sin^4\left(\theta+4 \sigma\right)\sin\left(\sigma-\theta\right)}{\sin^5\left(5\sigma\right)}=2\xi
\end{equation}
Let's  write:
\[f(\sigma)=g_1(\sigma)g_2(\sigma) \]
Where
\[ g_1(\sigma)=\frac{\sin^4\left(\theta+4 \sigma\right)}{\sin^4\left(5\sigma\right)} \qquad \textrm{and}\qquad g_2(\sigma)=\frac{\sin\left(\sigma-\theta\right)}{\sin\left(5\sigma\right)} \]
$g_1$ is positive and stricly increasing in $[-\frac{\theta}{4},0[$. Likewise $g_2$ is positive and stricly increasing since:
\[g_2'(\sigma)=\frac{\sin(4\sigma+\theta)+4\sin(\theta-\sigma)\cos(5\sigma)}{\sin^2\left(5\sigma\right)}>0\]
Therefore $f$ is continuous and strictly increasing. In addition:
\begin{equation}
f\left(-\frac{\theta}{4}\right)=0 \qquad  \textrm{ and }\qquad \lim_{\sigma\to 0^-}f\left(\sigma\right)=+\infty 
\end{equation}
Consequently, there exists a unique $\sigma \in [-\frac{\theta}{4},0[$ satisfying (15) or a unique root $y^*$ such as:
\begin{equation}
-\frac{\theta}{4} \leq Arg(y^*)<0
\end{equation}
\subsection{Absolute error theorem}
In this section, we will prove the absolute error theorem: $\forall\: \xi \in \R^{+*}$ and $\forall \:\theta \in \left[0,\pi/5\right]$: \[\left|y_1-y^{*}\right| <4.32 \textrm{ } 10^{-3}\]
We will start by discussing the unique properties of the selected starting point. 
\subsubsection{Starting point properties}

\begin{itemize}
\item P.1
\begin{equation}
|y^*| \leq |y_0|
\end{equation}
\item P.2  
\begin{equation}
|u+y^*|  \leq |u+y_0|
\end{equation}
\item P.3 Define
\[h_0=\sqrt[4]{\frac{u+y_0}{u+y^*}}\] 
\begin{equation}
\left|Arg\left(h_0\right) \right|\leq\frac{\pi}{20}   
\end{equation}
\end{itemize}
{\bf Proof:}
\begin{itemize}
\item {\bf P.1: ${\bf |y^*| \leq |y_0|}$} \\
\\
Using (3)
\[|{y^*}|^8|u+y^*|^2=4\xi^2\]
Therefore
\[|{y^*}|^8\left|1+|{y^*}|^2+2|{y^*}|\cos(\theta-\sigma)\right|=4\xi^2\]
\\
Since $1+|{y^*}|^2 \geq 2|{y^*}|$:
\[|{y^*}|^9\cos^2\left(\frac{\theta-\sigma}{2}\right)\leq \xi^2 \qquad \textrm{or}\qquad |{y^*}| \leq \left( \frac{\xi^2}{\cos^2\left(\frac{\theta-\sigma}{2}\right)}\right)^{\frac{1}{9}}\]
Also:
\[0 \leq \theta \leq \frac{\pi}{5} \qquad \textrm{and}\qquad -\frac{\theta}{4} \leq \sigma \leq 0\]
Which leads to: 
 \[0 \leq \frac{\theta-\sigma}{2} \leq \frac{\pi}{8} \] 
 and
 \[ \left|\cos\left(\frac{\theta-\sigma}{2}\right)\right| \geq \cos\left(\frac{\pi}{8}\right)=\sqrt{\frac{2+\sqrt{2}}{4}}=\alpha\]
It follows that:
\[|y^*| \leq \left(\frac{\xi}{\alpha}\right)^{\frac{2}{9}}=|y_0|\]
\item {\bf P.2: ${\bf |u+y^*|  \leq |u+y_0|}$} \\
\\
Since $|u|=1$ and the angle between $u$ and $y^*$ is  $\theta-\sigma$:
\[|u+y^*|^2=1+|y^*|^2+2|y^*|\cos\left(\theta-\sigma\right)\] 
Likewise, since the angle between $u$ and $y_0$ is  $\theta$:
\[|u+y_0|^2=1+|y_0|^2+2|y_0|\cos\left(\theta\right)\] 
From P.1: $|y^*|\leq|y_0|$. Also $0 \leq \cos\left(\theta\right)\leq \cos\left(\theta-\sigma\right)$. Therefore:
\[|u+y^*|  \leq |u+y_0|\]
\item {\bf P.3:} \\
\[{\bf h_0=\sqrt[4]{\frac{u+y_0}{u+y^*}} \qquad \left|Arg\left(h_0\right) \right|\leq\frac{\pi}{20}}\]
\\
Notice $u+y_0=(y_0+\cos(\theta))+ i \sin(\theta)$. Since $y_0>0$
\[0\leq \frac{\sin(\theta)}{y_0+\cos(\theta)} \leq \frac{\sin(\theta)}{\cos(\theta)}\]
Therefore
\[0\leq Arg\left({u+y_0}\right) \leq \theta \qquad \qquad \qquad  (P.3.1)\]
On the other hand, since ${y^*}^4(u+y^*)=2\xi \in \R^+$ and $Arg({y^*})=\sigma$:
\[Arg\left({u+y^*}\right)=-4\sigma\]
Since:
\[-\frac{\theta}{4}\leq \sigma \leq 0\]
\\
\[0\leq Arg\left({u+y^*}\right)\leq \theta \qquad \qquad \qquad (P.3.2)\]
From (P.3.1) and (P.3.2)
\[-\theta \leq Arg\left({u+y_0}\right)-Arg\left({u+y^*}\right)\leq \theta\]
Or
\[ \left|Arg\left(\frac{u+y_0}{u+y^*}\right) \right| \leq \theta \leq \frac{\pi}{5}\]
Since 
\[Arg(h_0)= \frac{1}{4} Arg\left(\frac{u+y_0}{u+y^*}\right) \]
\[\left|Arg\left(h_0\right) \right|\leq\frac{\pi}{20}\]
\end{itemize}
\subsubsection{Proof of absolute error theorem} 
\begin{equation}
{\bf \delta_1=|y_{1}-y^*|<C_0 \qquad \textrm{where } \qquad C_0=4.32 \: 10^{-3}}
\end {equation}
For $k \in \N$ define:
\begin{equation}
t_k=\sqrt[4]{\frac{2\xi}{u+y_k}}, \qquad h_k=\frac{y^*}{t_k}=\sqrt[4]{\frac{u+y_k}{u+y^*}}, \qquad \delta_k=|y_{k}-y^*|
\end {equation}
where $(y_k)_{k \in \N}$ is defined in (9) and (10) (section 3.1).\\
\\
Also for $t \in \C$ define:
\begin {equation}
P(t)=2 t^3+ \frac{2u}{5} t^2+\frac{u^2}{25}t
\end {equation}
Recall that ${y^*}=G\left(\xi,{y^*}\right)$:
\[{y^*}=\sqrt[5]{2\xi+\frac{2u^2}{5}{y^*}^3 +\frac{2u^3}{25}{y^*}^2+\frac{u^4}{125}{y^*}+\frac{u^5}{3125}}-\frac{u}{5} \]
Using (23), this can be simplified  to:
\[{y^*}+\frac{u}{5}=\sqrt[5]{2\xi+\frac{u^5}{3125}+\frac{u^2}{5}P(y^*)} \]
Or:
\begin {equation}
2\xi+\frac{u^5}{3125}=\left({y^*}+\frac{u}{5}\right)^5-\frac{u^2}{5}P(y^*)
\end {equation}
Likewise ${y_1}=G\left(\xi,{y_0}\right)$ leads to:
\[y_1+\frac{u}{5}=\sqrt[5]{2\xi+\frac{u^5}{3125}+\frac{u^2}{5}P(t_0)}= \sqrt[5]{\left({y^*}+\frac{u}{5}\right)^5+\frac{u^2}{5}\left(P(t_0)-P(y^*)\right)}\]
Therefore:
\begin {equation}
y_{1}+\frac{u}{5}=\left(y^*+\frac{u}{5}\right)\sqrt[5]{1+\epsilon_1} \qquad \textrm{ where } \qquad \epsilon_1=\frac{u^2}{5}\frac{P(t_0)-P(y^*)}{\left(y^*+\frac{u}{5}\right)^5}
\end {equation}
In the following two lemmas we will show that $|\epsilon_1|<0.049444787$
\\
\\
{\bf Lemma 1: Upper bound for $|\epsilon_1|$\\ \\
 Define:}
\[{\bf\omega=\frac{1}{2\alpha} \left(\left(\frac{y^*}{u}\right)^{\frac{1}{2}}+\left(\frac{u}{y^*}\right)^{\frac{1}{2}}\right) \qquad r=|y^*| \qquad d_1=3.82975138}\]
{\bf and}
\[
{\bf A=\frac{1}{20 d_1 \alpha^2} \qquad B=\left|\frac{u^\frac{1}{9}\omega^{\frac{2}{9}}-1}{\omega^2}\right|\qquad C=\frac{6 r^4+\frac{34}{5}r^3+\frac{21}{25}r^2+\frac{1}{25}r}{\left|y^*+\frac{u}{5}\right|^{5}}}
\]
${\bf  \epsilon_1}$ {\bf satisfies:}
\[{\bf |\epsilon_1|<ABC}\]
{\bf Proof:}\\
\begin {equation}
P(t_0)-P(y^*)=({t_0}-{y^*})\left(2{t_0}^2+2{t_0}{y^*}+2{y^*}^2+\frac{2}{5}({t_0}+{y^*})+\frac{u^2}{25}\right)
\end {equation}
From P.2, $|u+y^*|  \leq |u+y_0|$ therefore:  
\[|h_0|=\left|\frac{y^*}{t_0}\right|=\left|\sqrt[4]{\frac{u+y_0}{u+y^*}}\right| \geq 1\]
Or
\[|t_0| \leq r=|y^*|\]
This combined with (26) leads to:
\begin {equation}
\left|P(t_0)-P(y^*)\right| \leq \left|t_0-y^*\right|\left(6 r^2+\frac{4}{5}r+\frac{1}{25}\right)
\end {equation}
On the other hand:
\begin {equation}  
t_0-y^*=\frac{{t_0}^4-{y^*}^4}{{t_0}^3+{t_0}^2y^*+t_0{y^*}^2+{y^*}^3}
\end {equation}
Since 
\[{t_0}^4=\frac{2\xi}{u+y_0} \qquad \textrm{and}\qquad {y^*}^4=\frac{2\xi}{u+y^*}\] 
\[{t_0}^4-{y^*}^4=\frac{2\xi}{u+y^*}\frac{y^*-y_0}{u+y_0}={y^*}^4\frac{y^*-y_0}{u+y_0}\] 
It follows that:
\[t_0-y^*={y^*}^4\frac{y^*-y_0}{(u+y_0)({t_0}^3+{t_0}^2y^*+t_0{y^*}^2+{y^*}^3)}\]
Recall from (22):
\[t_0=\frac{y^*}{h_0}  \qquad \textrm{and} \qquad u+y_0={h_0}^4(u+y^*)\] 
Therefore:
\begin {equation}
t_0-y^*=y^*\frac{(y^*-y_0)}{(u+y^*)\left({h_0} +{h_0}^2 +{h_0}^3+{h_0}^4\right)}
\end {equation}
Since $ |h_0| \geq 1$, $|Arg(h_0)| \leq \pi/20$ and the modulus of a complex number is always greater than its real part \footnote{This is a good approximation due to the benefit of working near the real axis.}:
\begin {equation}
|{h_0}||1+{h_0}+{h_0}^2+{h_0}^3| \geq \sum_{k=0}^{3} \cos\left(k\frac{\pi}{20}\right) \geq d_1=3.82975138
\end {equation}
From (29):
\[|t_0-y^*|\leq|y^*|\frac{|y^*-y_0|}{d_1|u+y^*|}=\left|\frac{y_0}{y^*}-1\right|\frac{|y^*|^2}{d_1|u+y^*|}\]
Therefore
\begin {equation}
|t_0-y^*|\leq\left|\frac{y_0}{y^*}-1\right|\frac{r^2}{d_1|u+y^*|}
\end {equation}
Let's now turn our attention to the denominator of $\epsilon_1$ in equation (25). Since $\cos(u,y^*) \geq \cos(\pi/4)$:
\begin {equation}
\left|y^*+\frac{u}{5}\right| \geq \sqrt{r^2 +\frac{1}{25}+r \frac{\sqrt{2}}{5} }
\end {equation}
Using (25), (27) and (31):
\begin {equation}
|\epsilon_1| \leq \frac{1}{5 d_1}\frac{|\frac{y_0}{y^*}-1|}{|u+y^*|}\frac{6 r^4+\frac{4}{5}r^2+\frac{1}{25}r^2}{\left|y^*+\frac{u}{5}\right|^{5}}
\end {equation}
Going back to equation (3):
\[{y^*}^4(u+{y^*})=2\xi=2\alpha {y_0}^{\frac{9}{2}}\]
By factorizing:
\[u^{\frac{1}{2}}{y^*}^{\frac{9}{2}}\left(\left(\frac{u}{y^*}\right)^{\frac{1}{2}}+\left(\frac{y^*}{u}\right)^{\frac{1}{2}}\right)=2\alpha {y_0}^{\frac{9}{2}}\]
Since 
\[\omega=\frac{1}{2\alpha} \left(\left(\frac{y^*}{u}\right)^{\frac{1}{2}}+\left(\frac{u}{y^*}\right)^{\frac{1}{2}}\right) \]
\[2\alpha {y_0}^{\frac{9}{2}}=u^{\frac{1}{2}}{y^*}^{\frac{9}{2}}\left(\left(\frac{y^*}{u}\right)^{\frac{1}{2}}+\left(\frac{u}{y^*}\right)^{\frac{1}{2}}\right)=2\alpha u^{\frac{1}{2}}{y^*}^{\frac{9}{2}}\omega\]
Therefore
\[\left(\frac{y_0}{y^*}\right)^\frac{9}{2}=u^\frac{1}{2}\omega\]
Since $|Arg\left(\frac{y_0}{y^*}\right)|\leq \frac{\pi}{20}$
\begin {equation} 
\frac{y_0}{y^*}=u^\frac{1}{9}\omega^{\frac{2}{9}}
\end {equation}
From (33) and (34):
\[|\epsilon_1| \leq \frac{1}{5 d_1 }\frac{\left|u^\frac{1}{9}\omega^{\frac{2}{9}}-1\right|}{\left|(1+\frac{y^*}{u})\right|}\frac{6 r^4+\frac{4}{5}r^3+\frac{1}{25}r^2}{\left|y^*+\frac{u}{5}\right|^{5}}\]
Multiplying the numerator and the denominator by $|1+u/y^*|$
\begin{equation}
|\epsilon_1| \leq \frac{1}{5 d_1 }\frac{\left|u^\frac{1}{9}\omega^{\frac{2}{9}}-1\right|}{\left|(1+\frac{y^*}{u})(1+\frac{u}{y^*})\right|}\frac{6 r^4+\frac{4}{5}r^3+\frac{1}{25}r^2}{\left|y^*+\frac{u}{5}\right|^{\frac{5}{2}}}|1+\frac{u}{y^*}|
\end{equation}
Since 
\[\omega^2=\frac{1}{4\alpha^2}\left(1+\frac{y^*}{u}\right)\left(1+\frac{u}{y^*}\right)\]
\begin{equation}
|\epsilon_1| \leq \frac{1}{20 d_1 \alpha^2 }\left|\frac{u^\frac{1}{9}\omega^{\frac{2}{9}}-1}{\omega^2}\right|\frac{6 r^4+\frac{4}{5}r^3+\frac{1}{25}r^2}{\left|y^*+\frac{u}{5}\right|^{5}}|1+\frac{u}{y^*}|
\end{equation}
Recall that $|1+\frac{u}{y^*}|\leq (1+1/r)$, then:
\[\left(6 r^4+\frac{4}{5}r^3+\frac{1}{25}r^2\right)\left|1+\frac{u}{y^*}\right|\leq 6 r^4+\frac{34}{5}r^3+\frac{21}{25}r^2+\frac{1}{25}r\]
This combined with (36) leads to:
\begin{equation}
|\epsilon_1|\leq A B C 
\end{equation}
With 
\[ A=\frac{1}{20 d_1 \alpha^2 } \qquad B=\left|\frac{u^\frac{1}{9}\omega^{\frac{2}{9}}-1}{\omega^2}\right|\qquad C=\frac{6 r^4+\frac{34}{5}r^3+\frac{21}{25}r^2+\frac{1}{25}r}{\left|y^*+\frac{u}{5}\right|^{5}}\]
{\bf Lemma 2:  Upper bounds of A, B, C and ${\bf \epsilon_1}$}
\begin {itemize}
\item P.4
\[\bf A \leq 0.0152957\]
\item P.5
\[{\bf B \leq 0.1635792}\]
\item P.6
\[{\bf C\leq 19.7616704}\]
\item P.7
\[{\bf \epsilon_1 \leq 0.049444787}\]
\end{itemize}
\newpage
{\bf Proof:}\\
\begin {itemize}
\item P.4
\[\bf A \leq 0.0152957\]
Recall:
\[A=\frac{1}{20 d_1 \alpha^2}\]
Since: 
\[d_1=3.829751381 \qquad \alpha=\sqrt{\frac{2+\sqrt{2}}{4}}\]
\begin{equation}
{\bf A \leq 0.0152957}
\end{equation}
\item P.5
\[{\bf B \leq 0.1635792}\]
Define \[v=u^\frac{1}{9}\omega^{\frac{2}{9}}=y_0/y^*\] 
Using P.1  $|v| \geq 1$. Also since $y_0 \in \R^+$  \[\beta=arg(v)=-arg(y^*)=-\sigma\]
Therefore $0\leq \beta \leq \theta/4 \leq \frac{\pi}{20}$.
\begin{equation}
B=\left|\frac{v-1}{v^9}\right|=\sqrt{\frac{|v|^2+1-2 |v|\cos(\beta)}{|v|^{18}}}
\end{equation}
This expression can be rewritten as:
\begin{equation}
B=\sqrt{\frac{m^2}{(m+\cos(\beta))^{18}}+\frac{\sin^2(\beta)}{|v|^{18}}} \qquad \textrm{with}\qquad m=|v|-\cos(\beta) \geq 0
\end{equation}
Leading to:
\begin{equation}
B \leq \sqrt{(m^{\frac{8}{9}}+\cos(\beta)m^{-\frac{1}{9}})^{-18}+\sin^2(\beta)} 
\end{equation}
The maximum of the right-hand side is reached when the derivative with respect to $m$ is zero at $m=\cos(\beta)/8$. This leads to:
\begin{equation}
B \leq \sqrt{\frac{1}{9^2}\left(\frac{8}{9\cos(\beta)}\right)^{16}+\sin^2(\beta)}
\end{equation}
As a result:
\begin{equation}
{\bf B \leq 0.1635792}
\end{equation}
\item P.6
\[{\bf C\leq 19.7616704}\]
Recall that:
\[C=\frac{6 r^4+\frac{34}{5}r^3+\frac{21}{25}r^2+\frac{1}{25}r}{\left|y^*+\frac{u}{5}\right|^{5}}\]
Using (32)
\[C\leq \frac{6 r^4+\frac{34}{5}r^3+\frac{21}{25}r^2+\frac{1}{25}r}{\left(r^2+\frac{1}{25}+r\frac{\sqrt{2}}{5}\right)^{\frac{5}{2}}}\]
Therefore:
\begin{equation}
 C\leq 6 I_4+\frac{34}{5}I_3+\frac{21}{25}I_2+\frac{1}{25}I_1
\end{equation}
Where
 \[I_k=\frac{r^k}{\left(r^2 +\frac{1}{25}+r \frac{\sqrt{2}}{5}\right)^{\frac{5}{2}}} \qquad k=1,2,3,4\]
or
 \[I_k=\left(r^{\frac{10-2k}{5}} +\frac{1}{25}r^{-\frac{2k}{5}}+r^{\frac{5-2k}{5}} \frac{\sqrt{2}}{5}\right)^{-\frac{5}{2}} \]
By deriving with respect to $r$, the upper bound for $I_k$ is reached at:
\[r_k=\frac{(2k-5)\sqrt{2}+\sqrt{50-8k^2+40k}}{100-20k}\]
 Which leads to:
\begin{equation}
{\bf C\leq 19.7616704}
\end{equation}
\item P.7
\[{\bf \epsilon_1 \leq 0.049444787}\]
Using (37), (38), (43) and (45): 
\begin{equation}
\boxed {{\bf |\epsilon_1| \leq 0.049444787}}
\end{equation}
\end{itemize}
 {\bf Lemma 3: Upper bound of $|y_1-y^*|$\\\\
 Define:}
\[ {\bf d_2=4.801223053 \qquad C'=\frac{6 r^4+\frac{34}{5}r^3+\frac{21}{25}r^2+\frac{1}{25}r}{\left|y^*+\frac{u}{5}\right|^{4}} }\]
\begin{itemize}
\item P.8
\[{\bf \left|y_1-y^*\right|\leq \frac{1}{d_2}A.B.C'}\]
\item P.9
\[{\bf C' \leq 8.288300}\]
\item P.10
\[{\bf |y_1-{y^*}|\leq 0.004319271}\]
\end{itemize}
{\bf  Proof:}\\
\begin{itemize}
\item P.8
\[{\bf \left|y_1-y^*\right|\leq \frac{1}{d_2}A.B.C'}\]
From (25):
\begin{equation}
\left|y_1-y^*\right| =\left|\left(y^*+\frac{u}{5}\right)(\sqrt[5]{1+\epsilon_1}-1)\right|
\end{equation}
Using the maximum of the modulus of the derivative: 
\[|(\sqrt[5]{1+\epsilon_1}-1)| \leq \frac{1}{5\left(1-|\epsilon_1|\right)^\frac{4}{5}}|\epsilon_1|
 \]
Since $|\epsilon_1| \leq 0.049444787$:
\[\frac{1}{5\left(1-|\epsilon_1|\right)^\frac{4}{5}}\leq \frac{1}{d_2} \qquad  \textrm{where}\qquad  
d_2= 4.801223053\]
Then:
\begin{equation}
|(\sqrt[5]{1+\epsilon_1}-1)|  \leq \frac{1}{d_2}|\epsilon_1|
\end{equation}
and
\begin{equation}
\left|y_1-y^*\right| \leq \frac{1}{d_2}\left|y^*+\frac{u}{5}\right||\epsilon_1|
\end{equation}
Recall from Lemma 1 that $|\epsilon_1| \leq ABC$, therefore:
\begin{equation}
\left|y_1-y^*\right|\leq \frac{1}{d_2}A.B.C'
\end{equation}
With
\begin{equation}
C'=\left|y^*+\frac{u}{5}\right| C=\frac{6 r^4+\frac{34}{5}r^3+\frac{21}{25}r^2+\frac{1}{25}r}{\left|y^*+\frac{u}{5}\right|^4}
\end{equation}
\item P.9
\[{\bf C' \leq 8.288300}\]
Using (32)
\[C'\leq \frac{6 r^4+\frac{34}{5}r^3+\frac{21}{25}r^2+\frac{1}{25}r}{\left(r^2 +\frac{1}{25}+r \frac{\sqrt{2}}{5}\right)^2}\leq 6+a_3 J_3+ a_2 J_2+ a_1 J_1+ a_0 J_0\]
Where
\[a_3=\frac{34-12\sqrt{2}}{5} \qquad a_2=-\frac{3}{25} \qquad a_1=\frac{5-12\sqrt{2}}{125} \qquad a_0=-\frac{6}{625}\]
Where
\[J_k=\frac{r^k}{\left(r^2 +\frac{1}{25}+r \frac{\sqrt{2}}{5}\right)^2}=\left(r^{\frac{1}{2}} +\frac{1}{25}r^{-\frac{3}{2}}+r^{-\frac{1}{2}} \frac{\sqrt{2}}{5}\right)^{-2}\]
Since $a_0$, $a_1$ and $a_2$ are negative
\[C'\leq 6+a_3 J_3\]
Using the derivative, $J_k$'s maximum is reached at:
\[r_k=\frac{(k-2)\sqrt{2}+\sqrt{8+8k-2k^2}}{40-10k}\]
Therefore when $k=3$: $J_3 \leq 2.288300$ which leads to:
\begin{equation}
{\bf C' \leq 8.288300}
\end{equation}
\item P.10
\[{\bf |y_1-{y^*}|\leq 0.004319271}\]
From (37), (38), (43) and (52)
\begin{equation}
\boxed{{\bf |y_1-{y^*}|\leq 0.004319271}}
\end{equation}
\end{itemize}
\subsection{Relative error theorem}
{\bf Theorem 2:\\
\begin{equation}
{\bf \left|\frac{y_{1}}{y^*}-1\right|\leq C_1 \qquad \textrm{where } \qquad C_1=2.51 \: 10^{-2}}
\end {equation}}
{\bf Proof:}\\
Using (49):
\begin{equation}
\left|\frac{y_{1}}{y^*}-1\right|=\left|\frac{y_{1}-y^*}{y^*}\right| \leq \frac{1}{d_2}\left|\frac{\left(y^*+\frac{u}{5}\right)}{y^*} \epsilon_1\right|\end {equation}
Since $\epsilon_1 \leq ABC$:
\begin{equation}
\left|\frac{y_{1}}{y^*}-1\right| \leq \frac{1}{d_2}A B C''
\end {equation}
Where
\begin{equation} 
C''=C \frac{\left|y^*+\frac{u}{5}\right|}{r}=\frac{6 r^3+\frac{34}{5}r^2+\frac{21}{25}r+\frac{1}{25}}{\left|y^*+\frac{u}{5}\right|^4}
\end{equation}
Using (32)
\[C''\leq \frac{6 r^3+\frac{34}{5}r^2+\frac{21}{25}r+\frac{1}{25}}{\left(r^2+\frac{1}{25}+r \frac{\sqrt{2}}{5}\right)^2}=6 J_3+\frac{29}{5}J_2+\frac{21-5\sqrt{2}}{25}J_1+\frac{1}{\left(r^2 +\frac{1}{25}+r \frac{\sqrt{2}}{5}\right)}\]
 
\[\textrm{Since }\qquad \frac{1}{r^2 +\frac{1}{25}+r \frac{\sqrt{2}}{5}}\leq 25 \qquad \textrm{(decreasing function with maximum at 0)}\]
and $J_1 \leq 16.79664960$, $J_2 \leq 2.144660941$ and $J_3 \leq 0.671865984$ (using the maximum values established in P.9).
\begin{equation} 
C'' \leq  48.06772496
\end{equation}
 It follows from (56):
\begin{equation}
\boxed {{\bf \left|\frac{y_{1}}{y^*}-1\right|\leq 0.02504947}}
\end{equation}
\subsection{Speed of convergence}
{\bf Theorem 3: $\bf \forall k \geq 1$
\begin{equation}
\bf 
|y_{k+1}-y^*|<\frac{1}{K}|y_{k}-y^*| \qquad \textrm{and} \qquad |y_{k+1}-y^*|<\frac{C_0}{K^{k}}
\end {equation}}
{\bf with} $\bf  K=15.44$.
\\
\\
{\bf Proof:}\\
Recall that $r=|y^*|$ and from (22):
\[h_k=\frac{y^*}{t_k}=\sqrt[4]{\frac{u+y_k}{u+y^*}}\qquad \textrm{and} \qquad t_k=\sqrt[4]{\frac{2 \xi}{u+y_k}} \]
We will breakdown the proof into three lemmas:\\
\\
{\bf Lemma 4:}
For $k \in \N^*$ if $|y_k-y^*| \leq C_0$ then for $p=1,2,3,4$:
\[1-q \leq |{h_k}^p| \leq 1+q\]
\[|{t_k}^p|  \leq (1+q) r^p \]
with \[q=0.004338008\]
{\bf Proof:}\\
Since $\cos(u,y^*)\geq 0$ $|u+y^*| \geq |u|=1$ and 
\[\left|\frac{y_k-y^*}{u+y^*}\right|\leq |y_k-y^*|\] 
For $k\geq 1$ and $p=1,2,3,4$: using the maximum of the modulus of the derivative for the function: $x \longmapsto (1+x)^{\frac{p}{4}}$  
\[\left|{h_k}^{p}-1\right|=\left|\left(1+\frac{y_k-y^*}{u+y^*}\right)^{\frac{p}{4}}-1\right|\leq \frac{p}{4}\frac{|y_k-y^*|}{\left(1-|y_k-y^*|\right)^{1-\frac{p}{4}}}\leq\frac{|y_k-y^*|}{1-|y_k-y^*|} \]
Which leads to:
\[|{h_k}^p-1|\leq \frac{C_0}{1-C_0} \leq  q=0.004338008\]
Therefore
\begin{equation}
1-q \leq |{h_k}^p| \leq 1+q  
\end{equation}
Also, since $|{t_k}^p| = r^p|{h_k}^p|$
\[|{t_k}^p| \leq (1+q) r^p\]
In particular:
\[\left|{h_k} +{h_k}^2 +{h_k}^3+{h_k}^4\right|=\left|4+({h_k}-1) +({h_k}^2-1) +({h_k}^3-1)+({h_k}^4-1)\right|\]
Therefore
\[\left|{h_k} +{h_k}^2 +{h_k}^3+{h_k}^4\right|\geq 4-\left|{h_k}-1\right| -\left|{h_k}^2-1\right| -\left|{h_k}^3-1\right|-\left|{h_k}^4-1\right|\]
Using Lemma 4

\begin{equation}
\left|{h_k} +{h_k}^2 +{h_k}^3+{h_k}^4\right|\geq 4-4q
\end{equation}
\\
{\bf Lemma 5:}
\[{\bf \forall k\geq 1 \qquad \textrm{\bf If }\qquad  |y_k-y^*| \leq C_0 \qquad  \textrm{\bf then }\qquad |y_{k+1}-y^*|\leq \frac{|y_{k}-y^*|}{K}} \]
\\
{\bf Proof:}\\
Using exactly the same steps that led to (25):
\begin{equation}
y_{k+1}-y^*=\left(y^*+\frac{u}{5}\right)\left(\sqrt[5]{1+\epsilon_{k+1}}-1\right)
\end{equation}
With:
\[\epsilon_{k+1}=\frac{u^2}{5}\frac{P(t_k)-P(y^*)}{\left(y^*+\frac{u}{5}\right)^5}\]
Similarily, using the steps that led to (26) and (29), we obtain:
\begin{equation}
\epsilon_{k+1}=y^*\frac{u^2}{5}\frac{(y^*-y_k)(2{t_k}^2+2{t_k}{y^*}+2{y^*}^2+\frac{2}{5}\left({t_k}+{y^*}\right)+\frac{u^2}{25})}{(u+y^*)\left({h_k} +{h_k}^2 +{h_k}^3+{h_k}^4\right)\left(y^*+\frac{u}{5}\right)^{5}}
\end{equation}
Using Lemma 4: 
\[\left|\epsilon_{k+1}\right| \leq \frac{1+q}{5(4-4q)}\left(6I_3+\frac{4}{5}I_2+\frac{1}{25}I_1\right)C_0 \leq 0.000937238=q_1\]
Applying the maximum of the modulus of the derivative:
\begin{equation}
 |\sqrt[5]{1+\epsilon_{k+1}}-1| \leq \frac{1}{5(1- q_1)^{\frac{4}{5}}}|\epsilon_{k+1}|\leq \frac{1}{d_3}|\epsilon_{k+1}|
 \end{equation}
Where \[d_3=4.999128599\]
Using (63) and (65):
\[|y_{k+1}-y^*|\leq \frac{1}{d_3}\left|\left(y^*+\frac{u}{5}\right)\epsilon_{k+1}\right| \]
Or
\[|y_{k+1}-y^*|\leq \frac{1}{d_3}\left|y^*\frac{u^2}{5}\frac{(y^*-y_k)(2{t_k}^2+2{t_k}{y^*}+2{y^*}^2+\frac{2}{5}\left({t_k}+{y^*}\right)+\frac{u^2}{25})}{(u+y^*)\left({h_k} +{h_k}^2 +{h_k}^3+{h_k}^4\right)\left(y^*+\frac{u}{5}\right)^{4}}\right|\]
$|u+y^*|\geq 1$ and by using (32), (62) and Lemma 4:
\[|y_{k+1}-y^*|\leq |y_{k+1}-y^*|\frac{1+q}{20d_3(1-q) }\frac{6r^3+\frac{4}{5}r^2+\frac{1}{25}r}{\left(r^2+\frac{1}{25}+r\frac{\sqrt{2}}{5}\right)^{2}}\]
Therefore
\[|y_{k+1}-y^*|\leq |y_{k+1}-y^*|\frac{1+q}{20(1-q) d_3} \left(6J_3+\frac{4}{5}J_2+\frac{1}{25}J_1\right)\]
Since 
\[\frac{1+q}{20(1-q) d_3} \left(6J_3+\frac{4}{5}J_2+\frac{1}{25}J_1\right) \leq \frac{1}{K}\qquad  \textrm{where } \qquad K = 15.44198\]
\[
|y_{k+1}-y^*|\leq \frac{|y_{k}-y^*|}{K} 
\]
{\bf Lemma 6:}
$\forall k \geq 1$
\[|y_{k+1}-y^*|<\frac{1}{K}|y_{k}-y^*| \qquad \textrm{and} \qquad |y_{k+1}-y^*|<\frac{C_0}{K^{k}}\]
with $ K=15.44$.\\
\\
{\bf Proof:}\\
By induction, since  $|y_1-y^*|<C_0$ and $K>1$  and using Lemma 5: for all $k \geq 1$ $|y_k-y^*|<C_0$ and :
\[|y_{k+1}-y^*|\leq \frac{|y_{k}-y^*|}{K} \]
and
\[|y_{k+1}-y^*|\leq \frac{C_0}{K^{k}}\]

\section{Solving equations in Form 2}
Consider equations in Form 2: 
 \[\frac{z^5+z^4}{2}=\lambda\]
Recall from section 2 that a rotation $y=uz$ leads to Form 3: 
 \[\frac{y^5+uy^4}{2}=\xi\ \qquad \textrm{with}  \qquad \xi=|\lambda| \qquad \textrm{and} \qquad u =\left(\frac{|\lambda|}{\lambda}\right)^{\frac{1}{5}}\]
Therefore the same algorithm would apply to Form 2:
 \begin {enumerate}
 \item Starting point
\begin{equation} 
 z_0=\left(\frac{|\lambda|}{\alpha}\right)^\frac{2}{9}  \left(\frac{\lambda}{|\lambda|}\right)^{\frac{1}{5}} \qquad \textrm{with} \qquad \alpha=\sqrt{\frac{2+\sqrt{2}}{4}}  
\end{equation}
 \item Iteration 
\begin{equation} 
z_{k+1}=\sqrt[5]{2\lambda+\frac{2}{5}m_k^3 +\frac{2}{25}m_k^2+\frac{1}{125}m_k+\frac{1}{3125}}-\frac{1}{5}
\end{equation}
Where
\[m_k=\sqrt[4]{\frac{2\lambda}{1+z_k}}\]
\end {enumerate}
Since $|u|=1$ and $z_k=y_k/u$ for all $k \in \N$, the properties discussed in section (3) would apply:
\begin{enumerate}
\item Absolute error theorem: \[|z_{1}-z^*|<C_0 \qquad \textrm{where } \qquad C_0=4.32 \: 10^{-3}\]
\item Relative error theorem:
\[\left|\frac{z_1}{z^*}-1\right|<C_1 \qquad \textrm{where } \qquad C_1=2.51 \: 10^{-2}\]
\item Speed of convergence theorem: $\forall k \geq 1$
\[|z_{k+1}-z^*|<\frac{1}{K}|z_{k}-z^*| \qquad \textrm{and} \qquad |z_{k+1}-z^*|<\frac{C_0}{K^{k}}\]
with $ K=15.44$.
\\
\end{enumerate}
Also, since $Arg(z^*)=Arg(y^*)-\theta$, $|Arg(y^*)|\leq \pi/20$ and $|\theta|\leq \pi/5$:
\[\left|Arg(z^*)\right| \leq \frac{\pi}{4}\]

\section{Solving equations in Form 1}
\subsection{Challenge with iterative radicals}
Consider the equation
\[ x^5+x+a=0 \qquad a=0.01\]
One might think a trivial algorithm to find the real root of this equation is the following:
\[x_0=0 \qquad  \textrm{and} \qquad x_{k+1}=-\sqrt[5]{a+x}\]
Which produces the following sequence (approximated to the fifth digit):\\
$x_0=0$, $x_1=-0.3981$, $x_2=0.8275$, $x_3-0.9652$, $x_4=0.9909$, $x_5=-1.0002$, $x_6=0.9980$, $x_7=-1.0016$, $x_8=0.9983$, $x_9=-1.0017$, $x_{10}=0.9983$, $x_{11}=-1.0017$, $x_{12}=0.9983$, $x_{13}=-1.0017$, $x_{14}=0.9983$,...\\
\\
Notice the algorithm does not converge! But could another starting point $x_0$ lead to a different outcome? Outside the real-root itself, we will show below that it's not the case. \\
\\
Indeed assume that the sequence converges toward a real number $x^*$. Since \\$x^*({x^*}^4+1)=-a$, $|x^*| \leq a$.\\
\\
Use $\epsilon =0.01 a$. There exists an integer $k_0$ such that for all $k>k_0$: \[\left|x_{k}-x^*\right| <\epsilon\]
Note that: \[x^*=-\sqrt[5]{a+x^*}\]
\[x_{k+1}=-\sqrt[5]{a+x_{k}}\]
Using the mean value theorem for the function $-\sqrt[5]{a+x}$  between $x^*$ and $x_k$, there exist a real number $c$ between $x^*$ and $x_k$ such that:
\[\left|x_{k+1}-x^*\right|= \frac{1}{5}\left|a+c\right|^{-\frac{4}{5}}\left|x_{k}-x^*\right|\]
Since $\left|a+c\right|\leq|a|+|c|\leq| a|+|x^*|+\epsilon\leq 2a +\epsilon \leq 2.01 a$
\[\left|x_{k+1}-x^*\right|\geq 4.55 \left|x_{k}-x^*\right|\]
Which proves that the algorithm is not stable.
\subsection{Proposed algorithm}
Consider equation (1):
\[x^5+x+a=0 \qquad \textrm{with} \qquad a \neq 0\]
As discussed in section 2, a change of variable $z=a/x$ leads to: \[z^5+z^4=-a^4\]
As proven in section 4, the proposed algorithm provides an accurate first estimate $z_1$ for one of the five roots $z^*$:
\[|z_1-z^*|<C_0 \qquad \textrm{and}\qquad \left|\frac{z_1}{z^*}-1\right|<C_1\]
Consider $x^*=a/z^*$ and the sequence  $x_{k} =a/z_k$ for $k \in \N$ \\
\newpage
{\bf Lemma 7:}
\begin{enumerate}
\item Absolute error:
\[\left|x_1-x^*\right|<C_2 \qquad C_2=2.90  \: 10^{-2}\]
\item Relative error:
\[\left|\frac{x_1}{x^*}-1\right|<{C_1}' \qquad {C_1}'=2.57  \: 10^{-2}\]
\item Speed of convergence: $\forall k \geq 1$
\[\left|x_{k+1}-x^*\right|<\frac{1}{K'}\left|x_{k}-x^*\right| \qquad \textrm{and} \qquad \left|x_{k+1}-x^*\right|<\frac{C_2}{{K'}^k}\]
with $ K'=14.68$.
\end{enumerate}
{\bf Proof}:\\
\begin{enumerate}
\item
{\bf Absolute error:}
\[{\bf \left|x_1-x^*\right|<C_2 \qquad C_2=2.90  \: 10^{-2}}\]
\\
Notice:
\begin{equation}
\left|x_1-x^*\right|=\left|\frac{a}{z_1z^*}\right|\left|z_1-z^*\right|.
\end{equation}
\[{z^*}^4(1+z^*)=-a^4\]
Since, $\phi$ the argument $z^*$ is between $ [-\pi/4, \pi/4]$, 
\[\left|1+z^*\right|^2=1+|z^*|^2+2|z^*|\cos(\phi) \geq 1\]
Which leads to
\[\left|1+z^*\right| \geq 1\]
Since $|{z^*}|^4|1+z^*|=|a|^4$
\[|{z^*}|^4\leq|a|^4\]
Or
\begin{equation}
r=|z^*| \leq |a|.
\end{equation}
On the other hand, equation (2) implies:  
\begin{equation} 
|z^*|^5+|z^*|^4 \geq |a|^4.
\end{equation}
We distinguish two cases:
\begin{itemize}
\item $|a| \leq 1$: Using (69), $r\leq 1$. From (70): 
\[ |a|^4 \leq |z^*|^4+|z^*|^5 \leq 2 |z^*|^4 \qquad \textrm{therefore}\qquad \frac{|a|}{|z^*|} \leq 2^{\frac{1}{4}} \]
Since $|z_1/z^*-1|<C_1$ (relative error theorem in section 4)
\begin{equation}
\frac{1}{1+C_1}\leq |z^*/z_1|<\frac{1}{1-C_1}
\end{equation}
\[ \frac{|a|}{|z_1|} \leq \frac{|a|}{|z^*|} \frac{|z^*|}{|z_1|} \leq \frac{2^{\frac{1}{4}}}{1-C_1} \]
It follows that:
\[\left|x_1-x^*\right|=\left|\frac{a}{z_1z^*}\right|\left|z_1-z^*\right|=\left|\frac{a}{z_1}\right|\left|\frac{z_1}{z^*}-1\right|\]
Therefore:
\begin{equation}
\left|x_1-x^*\right| \leq \frac{2^{\frac{1}{4}}C_1}{1-C_1}\leq 2.90  \: 10^{-2}
\end{equation}

\item $|a|>1$: Using (70), $|z^*|>C_3=0.856$ (function $|z^*|^5+|z^*|^4$ is increasing in $|z^*|$ and ${C_3}^5+{C_3}^4 <1$\\
\\
Also from (70):
\[|a|^4 \leq |z^*|^5\left(1+\frac{1}{|z^*|}\right)\leq |z^*|^5\left(1+\frac{1}{C_3}\right)\]
Therefore
\begin{equation}
|z^*|  \geq |a|^\frac{4}{5} \sqrt[5]{\frac{C_3}{1+C_3}}
\end{equation}
Using (71):
\begin{equation}
|z_1|  \geq |a|^\frac{4}{5} \left(1-C_1\right)\sqrt[5]{\frac{C_3}{1+C_3}}
\end{equation}
From (73) and (74):
\[ \frac{|a|}{|z_1z^*|} \leq  |a|^{-\frac{3}{5}}\left(1-C_1\right)^{-1}\left(\frac{C_3}{1+C_3}\right)^{-\frac{2}{5}}<1.395955878\]
\begin{equation}
\left|x_1-x^*\right|\leq 0.005711884 \leq 2.90  \: 10^{-2}
\end{equation}
\end{itemize}
\item {\bf Relative error:}
\[\bf \left|\frac{x_1}{x^*}-1\right|<{C_1}' \qquad {C_1}'=2.57  \: 10^{-2}\]
Since $x^*=a/z^*$ and $x_1=a/x_1$:
\[ \left|\frac{x_1}{x^*}-1\right| = \left|\frac{z^*-z_1}{z_1}\right|=\left|\frac{z^*}{z_1}\right|\left|\frac{z_1}{z^*}-1\right|\]
Using the relative error theorem for Form 2 and (71):
\[\left|\frac{z_1}{z^*}-1\right|<C_1 \qquad \left|\frac{z^*-z_1}{z_1}\right|<\frac{1}{1-C_1}\]
Therefore
\[ \left|\frac{x_1}{x^*}-1\right| <\frac{C_1}{1-C_1}< {C_1}'=2.57  \: 10^{-2} \]
\item {\bf Speed of convergence:}
\[{\bf \left|x_{k+1}-x^*\right|<\frac{1}{K'}\left|x_{k}-x^*\right| \qquad \textrm{\bf with}\qquad  K'=14.68}\]
\\
Since $x=a/z$:
\[x_{k+1}-x^*=\frac{a}{z_{k+1}}-\frac{a}{z^*}=\frac{a(z^*-z_{k+1})}{z_{k+1}z^*}\] and
\[x_{k}-x^*=\frac{a}{z_{k}}-\frac{a}{z^*}=\frac{a(z^*-z_{k})}{z_{k}z^*}\]
Which leads to:
\[\left|\frac{x_{k+1}-x^*}{x_{k}-x^*}\right|=\left|\frac{z_{k}}{z^*}\right| \left|\frac{z^*}{z_{k+1}}\right|\left|\frac{z_{k+1}-z^*}{z_{k}-z^*}\right| \qquad \qquad (5.1)\]
First,
\[\left|\frac{z_{k+1}-z^*}{z_{k}-z^*}\right| \leq K\]
Second,
\[\left|\frac{z_{k}}{z^*}\right|=\left|1+\frac{z_{k}-z^*}{z^*}\right|\leq1+\frac{|z_{k}-z^*|}{|z^*|}\leq1+\frac{|z_{1}-z^*|}{|z^*|}\leq 1+C_1\]
Therefore:
\[\left|\frac{z_{k}}{z^*}\right|\leq 1+C_1 \qquad \qquad (5.2)\]
Likewise
\[\left|\frac{z_{k+1}}{z^*}\right|=\left|1+\frac{z_{k}-z^*}{z^*}\right| \geq1-\frac{|z_{k}-z^*|}{|z^*|}\geq1-\frac{|z_{1}-z^*|}{|z^*|}\geq 1-C_1\]
Therefore:
\[\left|\frac{z^*}{z_{k+1}}\right| \leq \frac{1}{1-C_1} \qquad \qquad (5.3)\]
From inequalities (5.1), (5.2) and (5.3)
\[\left|\frac{x_{k+1}-x^*}{x_{k}-x^*}\right|\leq \frac{1}{K}\frac{1+C_1}{1-C_1}\leq \frac{1}{K'}\qquad \textrm{with} \qquad K'=14.68
\]

\end{enumerate}

\section{Trigonometric algorithm}
The logic used in section 3.2 can be expanded to establish the geometric location of the five roots in the complex plane.\\
\\
\underline {Case 1: $0<\theta<\frac{\pi}{5}$}\\
\\
First, recall that ${y^*}=r e^{i \sigma}$  is a root of (3) if and only if two conditions are met:
\begin{equation}
f(\sigma)=\frac{\sin^4\left(\theta+4 \sigma\right)\sin\left(\sigma-\theta\right)}{\sin^5\left(5\sigma\right)}=2\xi  
\end{equation}
\begin{equation}
r=-\frac{\sin\left(\theta+4 \sigma\right)}{\sin\left(5\sigma\right)}>0
\end{equation}
Second, since \[{{y^*}}^5 \left(1+\frac{u}{{{y^*}}}\right)=2\xi\]
\[5 \sigma+Arg\left(1+\frac{u}{{y}}\right)=0 \mod 2\pi\]
Let $\sigma_k$ be the argument of the root ${{y^*}}_k$:
\[5 \sigma_k+Arg\left(1+\frac{u}{{y^*}_k}\right)=2k\pi \qquad -2\leq k \leq 2\]
When $\xi$ tends to $0^+$, either ${y^*}$ tends to ${y^*}_{-2}=-u$ with the argument $-\pi+\theta$ or ${y^*}$ tends to $0$ with the argument for that tends to: \[\frac{2k\pi}{5} \qquad \textrm{for} \qquad k=-1,0, 1,2\]
When $\xi$ tends $+\infty$, $|{y^*}_k|$ tends to $+\infty$ with an argument $\sigma_k$ that tends to: 
\[\frac{2k\pi-\theta}{4} \qquad \textrm{for} \qquad -2\leq k \leq 2\]
This intuition leads to considering the following intervals:  
\[I_{-1}=\left[-\frac{\pi}{2}-\frac{\theta}{4},-\frac{2\pi}{5}\right[  \qquad I_{0}=\left[-\frac{\theta}{4},-0\right[\]
 \[  I_{1}=\left[\frac{2\pi}{5}, \frac{\pi}{2}-\frac{\theta}{4}\right[  \qquad I_{2}=\left[\frac{4\pi}{5}, \pi-\frac{\theta}{4}\right[ \]
and 
\[I_{-2}=\left[-\pi +\theta,-\frac{4\pi}{5}\right[\]
$f$ is continuous in each of these intervals. $f$ is zero at the lower bound of $I_k$ for $k=-2,1,0$ and at the upper bound of  $I_k$ for $k=1,2$. $f$ also tends to $+\infty$ on the upper bound of  $I_k$ for $k=-2,1,0$  and at the lower bound $I_k$ for $k=1,2$. Notice also that if $\sigma \in I_k$ then (77) is satisfied ($r>0$).\\
\\
Consequently, equation (76) has at least one solution $\sigma_k$ in each of the five disjoint intervals. 
Since $f$ is strictly monotonic in each interval, we can define a bisection method to approximate $\sigma_k$ and the associated modulus from (77), leading to the five roots of equation (3). \\
\\
\underline {Case 2: $\theta=0$ or $\theta=\frac{\pi}{5}$}\\
\\
The same algorithm described in the previous case would apply to identify four roots, while the fifth can be obtained from Vieta's formulas. When $\theta=0$, we use a bisection in the intervals $I_{-2}$, $I_{-1}$, $I_{1}$ and $I_{2}$. For $\theta=\frac{\pi}{5}$, we use the bisection in the intervals  $I_{-1}$, $I_{0}$, $I_{1}$ and $I_{2}$.  

\section{Examples} 
All numbers presented in this section are approximated to the 10th decimal digit.
\subsection{Example 1: $\bf x^5+x+0.01=0 $}
\subsubsection{Trigonometric algorithm}
Using the transformation to Form (3)
\[\xi=0.000000005 \qquad \textrm{and}\qquad \theta=\frac{\pi}{5}\]
The bisection method in four intervals leads to roots ${y^*}_k$ for $k=-1,0,1,2$ of equation 3 and consequently to roots ${x^*}_k$ using the transformation ${x^*}_k=\frac{au}{{y^*}_k}$. The fifth root ($k=-2$) is obtained from Vieta's formula.\\
\\
\begin{tabular}{|R{0.5cm}|R{6cm}|R{6cm}|}
\hline $k$ & ${y^*}_k$ & $ {x^*}_k$   \\
\hline -2 & $-0.8090170025-0.5877852582i$ & $-0.0099999999$  \\
\hline -1 & $-0.0015494319-0.0098971415i$ & $-0.704595734+0.7071179873i$   \\
\hline 0 & $0.0098621565-0.0015443338i$ & $0.7095957339+0.7071176748i$   \\
\hline 1 & $0.0015788252+0.0098566936i$ & $ 0.7095957339-0.7071176748i$   \\
\hline 2 & $-0.0098915418+0.0015847876i$ & $ -0.704595734-0.7071179873i$   \\
\hline 
\end{tabular}
\\	
\begin{figure}[!h]
\centering
\begin{tikzpicture}
	\draw(-4,-3) rectangle (4,3);
     \draw[->] (-3, 0) -- (3, 0) node[right] {$\Re$};
  	\draw[->] (0, -2) -- (0, 2) node[above] {$\Im$};
  	\draw [fill=blue](-0.0099999999*2, 0) node[blue,  below left] {${x^*}_{-2}$}circle (2pt);
  	\draw [fill=blue](-0.704595734*2, 0.7071179873*2) node[blue,  below left] {${x^*}_{-1}$}circle (2pt);
	\draw [fill=blue](0.7095957339*2, 0.7071176748*2) node[blue,  below left] {${x^*}_{0}$}circle (2pt);
	\draw [fill=blue](0.7095957339*2, -0.7071176748*2) node[blue,  below left] {${x^*}_{1}$}circle (2pt);
	\draw [fill=blue](-0.704595734*2, -0.7071179873*2) node[blue,  below left] {${x^*}_{2}$}circle (2pt);
  \end{tikzpicture}
\caption{Roots of the equation $x^5+x+0.01=0$ }
\end{figure}
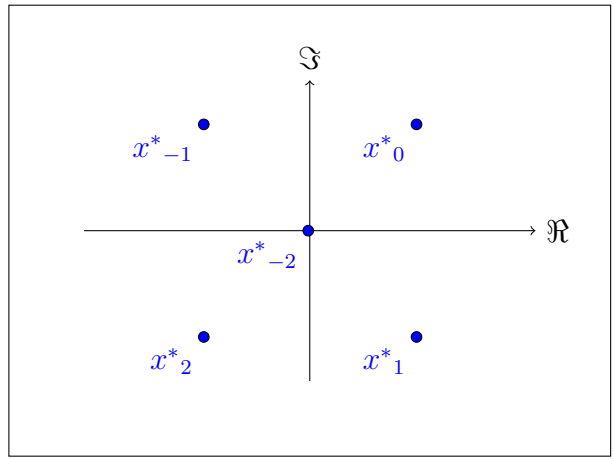
\subsubsection{Iteration of radicals algorithm}
Using the algorithm described in section 3 to estimate one root (coincides with $k=0$ from trigonometric algorithm): \[{y^*}_{0}=0.0098621565-0.0015443338i\]
\begin{tabular}{|R{1.5cm}|R{6cm}|R{2.4cm}|R{2.4cm}|}
\hline Iteration & $y_n$ & $|y_n-{y^*}_{0}|$ & $\left|\frac{y_n}{{y^*}_{0}}-1\right|$   \\
\hline 1 & $0.0098512048-0.0015389435i$ & $5.58 \times 10^{-4} $  &$2.18 \times 10^{-3}$\\
\hline 2 & $0.0098621666-0.0015443624i$ & $2.90 \times 10^{-6} $  & $1.13 \times 10^{-6}$ \\
\hline 3 & $0.0098621566-0.0015443337i$ & $1.51 \times 10^{-8}$  & $5.90 \times 10^{-9}$ \\
\hline 
\end{tabular}
\\	
\\
\\
Since $x=\frac{au}{y}$: ${x^*}_{0}=0.7095957339+0.7071176748i$ and the iterative estimates are:
\\
\\
\begin{tabular}{|R{1.5cm}|R{6cm}|R{2.4cm}|R{2.4cm}|}
\hline Iteration & $x_n$ & $|x_n-{x^*}_{0}|$ & $\left|\frac{x_n}{{x^*}_{0}}-1\right|$   \\
\hline 1 & $0.7106828395+0.707685341i$ & $1.23 \times 10^{-3} $  &$1.22 \times 10^{-3}$\\
\hline 2 & $0.7095928286+0.7071185567i$ & $3.04 \times 10^{-6} $  & $3.03 \times 10^{-6}$ \\
\hline 3 & $0.7095957376+0.7071176682i$ & $7.53 \times 10^{-9}$  & $7.51 \times 10^{-9}$ \\
\hline 
\end{tabular}
\\	
 
\subsection{Example 2: $\bf x^5+x+(3.08+1.68i)=0 $}
\subsubsection{Trigonometric algorithm}
Using the transformation to Form (3)
\[\xi \approx 75.75327872 \qquad \textrm{and}\qquad \theta \approx 0.228841153
\]
The bisection method in each of the intervals leads to roots ${y^*}_k$ of equation 3 and consequently to roots ${x^*}_k$ using the transformation ${x^*}_k=\frac{au}{{y^*}_k}$. \\
\\
\begin{tabular}{|R{0.5cm}|R{6cm}|R{6cm}|}
\hline $k$ & ${y^*}_k$ & ${x^*}_k$   \\
\hline -2 & $-2.4358363319-1.6419437613i$ & $-1.1834415151-0.1608289168i$  \\
\hline -1 & $0.6487113516-2.6125840601i$ & $-0.607389619+1.1531182439i$   \\
\hline 0 & $2.5580193297-0.0347499177i$ & $1.0110954185+0.9265109088i$   \\
\hline 1 & $0.6697215821+2.5335199517i$ & $1.116784747-0.7383651957i$   \\
\hline 2 & $-2.4145458637+1.5289087467i$ & $-0.3370490315-1.1804350402i$   \\
\hline
\end{tabular}
\\	
\begin{figure}[!h]
\centering
\begin{tikzpicture}
	\draw(-3.5,-2.6) rectangle (3.5,2.6);
     \draw[->] (-2.5, 0) -- (2.5, 0) node[right] {$\Re$};
  	\draw[->] (0, -2.0) -- (0, 2.0) node[above] {$\Im$};
  	\draw [fill=blue](-1.1834415151*1.5, -0.1608289168*1.5) node[blue,  below left] {${x^*}_{-2}$}circle (2pt);
  	\draw [fill=blue](-0.607389619*1.5, 1.1531182439*1.5) node[blue,  below left] {${x^*}_{-1}$}circle (2pt);
	\draw [fill=blue](1.0110954185*1.5, 0.9265109088*1.5) node[blue,  below left] {${x^*}_{0}$}circle (2pt);
	\draw [fill=blue](1.116784747*1.5, -0.7383651957*1.5) node[blue,  below left] {${x^*}_{1}$}circle (2pt);
	\draw [fill=blue](-0.3370490315*1.5, -1.1804350402*1.5) node[blue,  below left] {${x^*}_{2}$}circle (2pt);
  \end{tikzpicture}
\caption{Roots of the equation $x^5+x+(3.08+1.68i)=0$ }
\end{figure}
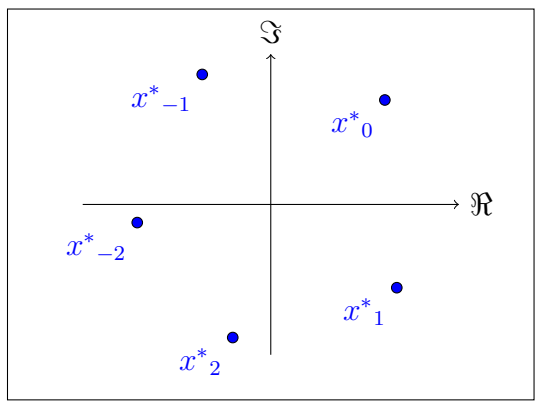
\subsubsection{Iteration of radicals algorithm}
Using the algorithm described in section 3 to estimate one root (coincides with $k=0$ from trigonometric algorithm): \[{y^*}_{0}=2.5580193297-0.0347499177i\]
\begin{tabular}{|R{2cm}|R{5.5cm}|R{2.4cm}|R{2.4cm}|}
\hline Iteration & $y_n$ & $|y_n-{y^*}_{0}|$ & $\left|\frac{y_n}{{y^*}_{0}}-1\right|$   \\
\hline 1 & $2.5575832547-0.0350982734i$ & $6.27 \times 10^{-4} $  &$2.13 \times 10^{-4}$\\
\hline 2 & $2.5580208152-0.0347474236i$ & $2.71 \times 10^{-6} $  & $9.23 \times 10^{-7}$ \\
\hline 3 & $2.5580193271-0.0347499325i$ & $1.18 \times 10^{-8}$  & $4.00 \times 10^{-9}$ \\
\hline 
\end{tabular}
\\	
\\
\\
Since $x=\frac{au}{y}$, ${x^*}_{0}=1.0110954185+0.9265109088i
$ and the iterative estimates are:
\\
\\
\begin{tabular}{|R{1.5cm}|R{6cm}|R{2.4cm}|R{2.4cm}|}
\hline Iteration & $x_n$ & $|x_n-{x^*}_{0}|$ & $\left|\frac{x_n}{{x^*}_{0}}-1\right|$   \\
\hline 1 & $1.0111375519+0.926807176i$ & $2.99 \times 10^{-4} $  &$4.07 \times 10^{-4}$\\
\hline 2 & $1.0110957554+0.9265093895i$ & $1.56 \times 10^{-6} $  & $2.12 \times 10^{-6}$ \\
\hline 3 & $1.0110954141+0.9265109156i$ & $8.09 \times 10^{-9}$  & $1.01 \times 10^{-8}$ \\
\hline 
\end{tabular}
\\	

\section{Future direction of work}
It would be interesting to study the feasibility of expanding the proposed iteration of radicals algorithm to identify all roots by selecting appropriate  starting points and the right combination of the 5th and 4th unity roots .


\begin{thebibliography}{9}

\bibitem[1]{1} Abel, Niels Henrik (1881) [1828], "Sur la resolution algébrique des équations", in Sylow, Ludwig; Lie, Sophus (eds.), Œuvres Complètes de Niels Henrik Abel (in French), II (2nd ed.), Grøndahl \& Søn, pp. 217–243
\bibitem[2]{2} Edwards, Harold M. (1984). Galois Theory. Springer-Verlag. ISBN 0-387-90980-X
\bibitem[3]{3} Al-Khwarizmi, Muhammad ibn Musa. (circa 825). Al-kitab al-mukhtasar fi hisab al-gabr wa’lmuqabala
\bibitem[4]{4} Cardano, Gerolamo (1545). Artis Magnæ
\bibitem[5]{5} O'Connor, John J.; Robertson, Edmund F., "Lodovico Ferrari", MacTutor History of Mathematics archive, University of St Andrews.
\bibitem[6]{6} Adamchik, Victor (2003). "Polynomial Transformations of Tschirnhaus, Bring, and Jerrard". ACM SIGSAM Bulletin. 
\bibitem[7]{7} Hermite, Charles (1858). "Sur la résolution de l'équation du cinquème degré". Comptes Rendus de l'Académie des Sciences. XLVI (I): 508–515.
\bibitem[8]{8} Brioschi, Francesco (1858). "Sul Metodo di Kronecker per la Risoluzione delle Equazioni di Quinto Grado". Atti Dell'i. R. Istituto Lombardo di Scienze, Lettere ed Arti. I: 275–282.
\bibitem[9]{9} Kronecker, Leopold (1858). "Sur la résolution de l'equation du cinquième degré, extrait d'une lettre adressée à M. Hermite". Comptes Rendus de l'Académie des Sciences. XLVI (I): 1150–1152.
\bibitem[10]{10} Cockle, James (1860). "Sketch of a Theory of Transcendental Roots". The London, Edinburgh, and Dublin Philosophical Magazine and Journal of Science.
\bibitem[11]{11} Harley, Robert (1862). "On the Transcendental Solution of Algebraic Equations". Quart. J. Pure Appl. Math. 5: 337–361
\bibitem[12]{12} Klein, Felix (1888). Lectures on the Icosahedron and the Solution of Equations of the Fifth Degree. Trübner \& Co. ISBN 978-0-486-49528-6.
\bibitem[13]{13} Doyle, Peter; Curt McMullen (1989). "Solving the quintic by iteration". Acta Math. 163: 151–180. 
\bibitem[14]{14} M.L.Glasser, (2000) “Hypergeometric functions and the trinomial equation”, Journal of Computational and Applied Mathematics, Volume 118, Issues 1-2, Pages 169-173


\end{thebibliography}
\end{document}